\documentclass[letterpaper, 10 pt, conference]{ieeeconf}  

\IEEEoverridecommandlockouts                              

\usepackage{afterpage}
\usepackage{epsfig} 
\usepackage{times} 
\usepackage{amsmath} 
\usepackage{amssymb}  
\usepackage{amsfonts}
\usepackage{mathptmx} 
\usepackage{tabularx}
\usepackage{epstopdf}
\usepackage{subfigure}
\usepackage{float}
\usepackage{epic,color}
\usepackage{mathrsfs}
\usepackage{dsfont,MnSymbol}
\usepackage{algorithm}
\usepackage{algorithmic}
\usepackage{multirow} 
\usepackage[T1]{fontenc}
\newcounter{rmnum}

\newcounter{anum}

\newlength{\noteWidth}
\long\def\notes#1{\ifinner
             {\tiny #1}
             \else
              \marginpar{\parbox[t]{\noteWidth}{\raggedright\tiny #1}}
               \fi}

\def\IEEEQEDclosed{\mbox{\rule[0pt]{1.3ex}{1.3ex}}}

\def\qed{\ifmmode\IEEEQEDclosed\else{\unskip\nobreak\hfil
\penalty50\hskip1em\null\nobreak\hfil\IEEEQEDclosed
\parfillskip=0pt\finalhyphendemerits=0\endgraf}\fi}

\def\qed{\hspace*{\fill}~\IEEEQED\par\endtrivlist\unskip}

\def\Re{\mathbb{R}}

\def\Theorem#1{Thm.~\ref{#1}}

\def\Sec#1{Sec.~\ref{#1}}

\def\notes#1{\marginpar{\tiny #1}\typeout{Notes!
Notes!
Notes!
}}
\renewcommand{\notes}[1]{\typeout{notes!}}

\def\Re{\field{R}}
\def\k{{\sf K}}

\def\Sec#1{Sec.~\ref{#1}}

\def\clZ{{\cal Z}}

\def\Sec#1{Sec~\ref{#1}}

\def\E{{\sf E}}
\def\Expect{{\sf E}}

\def\Expect{{\sf E}}

\def\Sec#1{Sec.~\ref{#1}}

\def\P{{\sf P}}

\def\IEEEQEDclosed{\mbox{\rule[0pt]{1.3ex}{1.3ex}}}
\def\qed{\nobreak\hfill\IEEEQEDclosed}

\def\clZ{{\cal Z}}

\newtheorem{theorem}{Theorem}

\newtheorem{remark}{Remark}
\newtheorem{proposition}{Proposition}

\def\beq{\begin{eqnarray}} 
\def\bc{\begin{center}} 
\def\be{\begin{enumerate}}
\def\bi{\begin{itemize}} 
\def\bs{\begin{small}}
\def\bS{\begin{slide}}
\def\ec{\end{center}} 
\def\ee{\end{enumerate}}
\def\ei{\end{itemize}}
\def\es{\end{small}}
\def\eS{\end{slide}}
\def\eeq{\end{eqnarray}}

\newcommand{\newP}[1]{\medskip\noindent{\bf #1:}}

\newcommand{\PP}{{\sf P}}

\newcommand{\ud}{\,\mathrm{d}}

\def\Re{\mathbb{R}}
\def\E{{\sf E}}

\def\Sec#1{Sec.~\ref{#1}}

\def\Prop#1{Prop.~\ref{#1}}

\def\Expect{{\sf E}}

\def\clZ{{\cal Z}}

\renewcommand{\Re}{\mathbb{R}}

\newcommand{\X}{X}

\title{\LARGE \bf
Derivation and Extensions of the Linear Feedback Particle Filter \\ based
on Duality Formalisms}

\author{Jin-Won Kim, Amirhossein Taghvaei and Prashant G. Mehta 
\thanks{Financial support from the NSF CMMI grant 1462773 is gratefully acknowledged. 
}
\thanks{J-W. Kim, A.~Taghvaei and P.~G.~Mehta are with the Coordinated
  Science Laboratory and the Department of Mechanical Science and
  Engineering at the University of Illinois at Urbana-Champaign
  (UIUC); Corresponding email: mehtapg@illinois.edu}
}

\begin{document}

\maketitle
\thispagestyle{empty}
\pagestyle{empty}

\begin{abstract}

This paper is concerned with a duality-based approach to derive the linear feedback
particle filter (FPF).  The FPF is a controlled
interacting particle system where the control law is designed to
provide an {\em exact} solution for 
the nonlinear filtering problem.  For the linear Gaussian special case,
certain simplifications arise whereby the linear FPF is identical to
the square-root form of the ensemble Kalman filter.  For this and for the
more general nonlinear non-Gaussian case, it has been an open problem
to derive/interpret the FPF control law as a solution of an optimal
control problem.  In this paper, certain duality-based arguments are
employed to transform the filtering problem into an optimal control
problem.  Its solution is shown to yield the deterministic form of the
linear FPF.  An extension is described to incorporate stochastic
effects due to noise leading to a novel homotopy of exact ensemble
Kalman filters.  All the derivations are based on duality formalisms.

\end{abstract}
\section{Introduction}

Amongst the many derivations of the Kalman filter, a particularly
appealing one (to control theorists) is based on duality
arguments; cf.,~\cite[Chapter 7]{astrom1970}.  In this derivation, the
linear Gaussian filtering (estimation) problem is modeled as an
optimization problem involving minimization of the mean-squared error
(variance) of the estimate.  By designing a suitable dual process (evolving in backward time), the
optimization problem is transformed into a finite-horizon
deterministic linear quadratic (LQ) optimal control problem.  The
solution to this LQ problem yields the Kalman filter.

The duality theorem originally appeared in Kalman's celebrated
paper~\cite{kalman1960}.  The duality relationship between optimal
control and linear filtering has proved to be very useful:
(i) it provides an interpretation of the Kalman
filter as the minimum variance estimator; (ii) it helps explain why the
filtering equation for the error covariance is same as the dynamic Ricatti
equation (DRE) of the optimal control; and 
(iii) quantitative results on the asymptotic properties of the
solution of the DRE have been used to derive results on
asymptotic stability of the linear filter~\cite{ocone1996}.

Given the
historical significance of these problems, several extensions have
been considered.  Early on, the
linear duality was related to the duality in nonlinear
programming~\cite{simon1970duality}. Extensions of the duality principle to
constrained linear estimation are described in~\cite{goodwin2005}.

Extensions of the duality theory to nonlinear
non-Gaussian settings is a closely related subject of historical importance beginning
with Mortensen's maximum-likelihood estimator~\cite{mortensen1968}.  Over the years,
there have been a number of important contributions in this area, 
e.g., the use of logarithmic (Hopf-Cole) transformation to convert the
filtering equation into the Hamilton-Jacobi-Bellman (HJB) equation of
optimal control~\cite{FlemingMitter82}.  An information theoretic interpretation
for the same appears in~\cite{mitter2003}.  More recently, these
considerations have led to the development and
application of the path integral approaches to filtering and smoothing
problems~\cite{ruiz_kappen2017}.  For a recent review with a
comprehensive reference list, cf.,~\cite{chetrite2015}.  

The goal of this paper is to generalize and apply Kalman's
duality principle for the purposes of deriving particle filters, in
the tractable linear Gaussian settings.  Specifically, the linear feedback particle
filter (also referred to as the {\em square-root} form of the ensemble
Kalman filter) is derived by extending results from the classical
duality theory.           

As some pertinent background, the feedback particle filter (FPF) is an
example of a controlled interacting particle system to approximate the
solution of the continuous-time nonlinear filtering problem.  In FPF,
the importance sampling step of the conventional particle filter is
replaced with feedback control. Other steps
such as resampling, reproduction, death or birth of particles are
altogether avoided. A salient feature of the FPF is that it is
an exact filter even in nonlinear-non-Gaussian settings.   An
expository review of the continuous-time filters including the
progression from the Kalman filter (1960s) to the ensemble Kalman
filter (1990s) to the feedback particle filter (2010s) appears
in~\cite{TaghvaeiASME2017}.

The contributions of this paper are as follows: This paper is the
first to present a derivation of linear FPF/ensemble Kalman filters
based on duality considerations.  The adaptation of the duality
formalisms to the particle filter, as proposed here, is also
original.  These considerations are used to obtain a novel homotopy
of stochastic and deterministic (linear) particle filters all of which are
exact. 

Although it is not an explicit focus of this paper, the derivation potentially
provides new tools to analyze and interpret the optimality properties of the filter, in terms of the simulation
variance of its estimates.  These properties are important in
applications of ensemble Kalman filters, in particular, for
high-dimension problems.

The outline of the remainder of this paper is as follows:
\Sec{sec:prelim} provides the background on FPF.  In
\Sec{sec:duality}, the duality formalism for the particle filter
is introduced to construct an optimization problem.  The solution to this
problem leading to the deterministic form of the linear FPF is
presented in \Sec{sec:det}.  The generalization to the stochastic case
appears on \Sec{sec:stoch}.     
All the proofs appear in the Appendix.

\section{Preliminaries and background}
\label{sec:prelim}

In this paper, we consider the linear Gaussian filtering problem:
\begin{subequations}
	\begin{align}
	\ud X_t &=  A_tX_t \ud t + \ud B_t\label{eq:dyn}\\
	\ud Z_t &= C_tX_t\ud t + \ud W_t\label{eq:obs}
	\end{align}
\end{subequations}
where $X_t \in \Re^d$ is the (hidden) state at time $t$, $Z_t \in \Re^m$ is the
observation; $A_t$, $C_t$ are matrices of appropriate
dimension whose elements are continuous in $t$; and $\{B_t\}$, $\{W_t\}$ are mutually independent Wiener
processes taking values in $\Re^d$ and $\Re^m$, respectively. The covariance associated with $\{B_t\}$
and $\{W_t\}$ are denoted by $Q_t$ and $R_t$, respectively.  
The initial condition $X_0$ is drawn from a Gaussian
distribution $\mathcal{N}(m_0,\Sigma_0)$, independent of $\{B_t\}$ or
$\{W_t\}$. It is assumed that the covariance matrices $Q_t$, $R_t$ and
$\Sigma_0$ are strictly positive definite for all $t$.  The filtering problem is to compute the posterior
distribution $\P(X_t|\clZ_t)$ where $\clZ_t:=\sigma(Z_s;s\in[0,t])$
denotes the filtration (the time-history of observations). 

For the linear Gaussian problem~\eqref{eq:dyn}-\eqref{eq:obs}, the posterior distribution
$\PP(\X_t|\clZ_t)$ is Gaussian $\mathcal{N}(m_t,\Sigma_t)$, whose mean
and covariance are given by the Kalman-Bucy filter~\cite{kalman1961}:
\begin{subequations}
\begin{align}
\ud m_t &= A_t m_t\ud t + \mathsf{K}_t(\ud Z_t - C_t m_t\ud t)\label{eq:KF-mean}\\
\frac{\ud \Sigma_t}{\ud t} &= A_t \Sigma_t +  \Sigma_tA_t^\top + Q_t - \Sigma_tC^\top R_t^{-1} C\Sigma_t\label{eq:KF-variance}
\end{align}
\end{subequations}
where $\k_t:=\Sigma_tC^\top R_t^{-1}$ is the Kalman gain and the filter is
initialized with the Gaussian prior $(m_0,\Sigma_0)$. 

Feedback particle filter (FPF) is a controlled interacting
particle system to approximate the Kalman filter\footnote{Although the
considerations of this paper are limited to the linear Gaussian
problem~\eqref{eq:dyn}-\eqref{eq:obs}, the FPF algorithm is more
broadly applicable to nonlinear non-Gaussian filtering
problems~\cite{taoyang_TAC12,yang2016}.}.  
In the following, the 
McKean-Vlasov stochastic differential equation (sde) are presented for the linear FPF algorithm.  For these models,
the state at time $t$ is denoted as $\bar{X}_t$. Two types of FPF
algorithm have been described in the literature~\cite{yang2016,AmirACC2016}: 

\newP{ (A) Stochastic linear FPF}  The state $\bar{X}_t$ evolves
according to the McKean-Vlasov sde:
\begin{align}
\ud \bar{X}_t &= A_t \bar{X}_t\ud t + \ud\bar{B}_t + \bar{\k}_t(\ud Z_t - \frac{C_t\bar{X}_t + C_t\bar{m}_t}{2}\ud t) \label{eq:Xbart-stochastic} 
\end{align}
where $\bar{\k}_t =\bar{\Sigma}_tC_t^\top R_t^{-1}$ is the Kalman gain; the
mean-field terms are the mean $\bar{m}_t=\Expect[\bar{X}_t|\clZ_t]$
and the covariance $\bar{\Sigma}_t
=\Expect[(\bar{X}_t-\bar{m}_t)(\bar{X}_t-\bar{m}_t)^\top|\clZ_t] $; 
$\{\bar{B}_t\}$ is an independent copy of the process noise $\{B_t\}$;
and the initial condition $\bar{X}_0 \sim \mathcal{N}(m_0,\Sigma_0)$.  

\newP{(B) Deterministic linear FPF} The McKean-Vlasov sde is:
\begin{align}
\ud \bar{X}_t =& A_t \bar{X}_t\ud t + \frac{1}{2}Q_t
\bar{\Sigma}_t^{-1}(\bar{X}_t-\bar{m}_t)\ud t + \bar{\k}_t(\ud Z_t - \frac{C_t\bar{X}_t + C_t\bar{m}_t}{2}\ud t)\label{eq:Xbart}
\end{align}
where (as before) $\bar{\k}_t =\bar{\Sigma}_tC_t^\top R_t^{-1}$ is the Kalman
gain; the mean $\bar{m}_t=\Expect[\bar{X}_t|\clZ_t]$
and the covariance $\bar{\Sigma}_t
=\Expect[(\bar{X}_t-\bar{m}_t)(\bar{X}_t-\bar{m}_t)^\top|\clZ_t]$; the initial condition $\bar{X}_0 \sim \mathcal{N}(m_0,\Sigma_0)$.

The difference between the two algorithms is that the process noise term $\ud \bar{B}_t$ in the stochastic FPF~\eqref{eq:Xbart-stochastic} is replaced by $\frac{1}{2}Q_t \bar{\Sigma}_t^{-1}(\bar{X}_t-\bar{m}_t)$ in the deterministic FPF~\eqref{eq:Xbart}.

The following Proposition, borrowed from~\cite{AmirACC2016}, shows
that both the filters are exact:

\medskip

\begin{proposition}
\label{thm_lin} {\em (Theorem~1 in~\cite{AmirACC2016})} Consider the linear Gaussian filtering problem~\eqref{eq:dyn}-\eqref{eq:obs}, and the linear
FPF (Eq.~\eqref{eq:Xbart-stochastic} or Eq.~\eqref{eq:Xbart}).  If $\PP(X_0)=\PP(\bar{X}_0)$ then
\begin{equation*}
\PP(\bar{X}_t|\clZ_t)=\PP(\X_t | \clZ_t),\quad \forall t>0
\end{equation*}
Therefore, $m_t = \bar{m}_t$ and $\Sigma_t = \bar{\Sigma}_t$. 
\label{thm:consistency-FPF-lin}
\end{proposition}

\medskip

In a numerical
implementation, the filter is simulated with $N$ interacting
particles where $N$ is typically large.  The filter state is $\{X_t^i:1\le i\le N\}$, where
$X_t^i$ is the state of the $i^{\text{th}}$-particle at time $t$.  The evolution of $X_t^i$ is
obtained upon empirically approximating the mean-field terms.  For
example, the finite-$N$ algorithm for the stochastic FPF is as follows:

\newP{(A) Finite-$N$ stochastic FPF} The evolution of $X_t^i$ is given by
the sde:
\begin{align}
\ud X^i_t &= A_t X^i_t\ud t + \ud B_t^i + \k^{(N)}_t(\ud Z_t -
\frac{C_tX^i_t + C_tm^{(N)}_t}{2}\ud t) 
\label{eq:Xit-s}
\end{align}
where $\k^{(N)}_t :=\Sigma^{(N)}_tC_t^\top R_t^{-1}$; $\{B^i_t\}_{i=1}^N$ are
independent copies of $B_t$; $X^i_0
\stackrel {\text{i.i.d}}{\sim} \mathcal{N}(m_0,\Sigma_0)$ for $i=1,2,\ldots,N$; and the
empirical approximations of the two mean-field terms are as follows:
\begin{align*}
m^{(N)}_t&:=\frac{1}{N}\sum_{j=1}^N X^i_t,~ \Sigma^{(N)}_t
:=\frac{1}{N-1}\sum_{j=1}^N (X^i_t-m^{(N)}_t)(X^i_t-m^{(N)}_t)^\top
\label{eq:empr_app_mean_var}
\end{align*}

The McKean-Vlasov sde~\eqref{eq:Xbart-stochastic} represents the
mean-field limit (as $N\rightarrow\infty$) of the finite-$N$
system~\eqref{eq:Xit-s}.

\section{Duality}
\label{sec:duality}

\subsection{Classical duality} 
Classical duality is concerned with the problem of constructing a
minimum variance estimator for the random variable $a^\top X_T$, where $X_T$ is the hidden state (defined
according to the model~\eqref{eq:dyn}) at (some fixed) time $T>0$ and $a\in\Re^d$ is an
arbitrary (but known) vector. 

Given observations $Z_t$ for $t\in[0,T]$ (defined according to the
model~\eqref{eq:obs}), the assumed linear structure for the causal estimator is as follows:
\begin{equation}
S_T = b_T^\top m_0 + \int_0^T u_t^\top \ud Z_t
\label{eq:KF_LP_est}
\end{equation}
The filter is thus parametrized by the vector $b_T\in
\Re^{d}$ and
the (possibly) time-dependent control input $u:[0,T]\to\Re^m$.  The
latter is denoted as $u_t$.  

The optimal filter parameters are obtained by solving the following
mean-squared optimization problem:
\[
\min_{b_T,u_t} \;\; \E [|S_T - a^\top X_T|^2]
\]  
subject to the following constraints:~\eqref{eq:dyn}-\eqref{eq:obs} for the
processes $X_t$ and $Z_t$, and~\eqref{eq:KF_LP_est} for the random
variable $S_T$.  

In the duality based derivation of Kalman filter, the optimization
problem is converted into a deterministic optimal control problem by
introducing a suitable dual process; cf.,~\cite[Chapter 7]{astrom1970}.  The solution of
the optimal control problem yields the optimal control input $u_t$ and
vector $b_T$.  By expressing the 
estimator~\eqref{eq:KF_LP_est} in its differential form, one obtains
the Kalman filter.

\subsection{Duality for the particle filter}

A particle filter empirically approximates the posterior
distribution.  By the $L^2$ optimality property of the posterior, the particle filter can be used to obtain a mean-squared
estimate of {\em any} arbitrary function of the hidden state.   

The particle analogue of~\eqref{eq:KF_LP_est} is as follows: 
Define $N$ random variables
$\{S_T^i:1\le i\le N\}$ according to
\begin{equation}
S_T^i = \sum_{j=1}^N(b_T^{ij})^\top X_0^j + \int_0^T (u_t^i)^\top \ud Z_t
\label{eq:FPF_LP_est_1}
\end{equation}
where the filter parameters now are $b_T^{ij}\in\Re^d$ for $i,j=1,\ldots,N$ and
the control input $u^i:[0,T]\to\Re^m$ for $i=1,\ldots,N$.  As
before, the filter state depends linearly on the known data -- observations up to time $T$ and the initial condition $X_0^j$ for
$j=1,\ldots,N$. The initial condition $X^j_0
\stackrel {\text{i.i.d}}{\sim} \mathcal{N}(m_0,\Sigma_0)$.  Such a sampling ensures that the filter provides
consistent estimates (as $N\rightarrow\infty$) at time $T=0$.

The optimal filter parameters are chosen
by solving the following mean-squared optimization problem:
\[
\min_{b_T^{ij},u_t^i} \;\; \E \Big[\big|\frac{1}{N}\sum_{i=1}^N f(S_T^i) -
f(a^\top X_T)\big|^2\Big]
\]  
subject to the constraints:~\eqref{eq:dyn}-\eqref{eq:obs} for the
processes $X_t$ and $Z_t$, and~\eqref{eq:FPF_LP_est_1} for the random
variables $S_T^i$. This is a multi-objective optimization
problem (as the function $f$ is arbitrary).  With the choice of $f(x)=x$, one obtains
the optimization problem considered in classical duality.  


There is a modeling trade-off here:  Ideally, one would want to
consider a large enough class of functions $f$, and the filter parameters
$b$ and $u$, that can represent and approximate the posterior
distribution.  However, the class may not be too large because one
would still want to be able to solve the multi-objective optimization problem.

\medskip

In order to make the analysis tractable, the following simplifications
are proposed:

\newP{Step 1} Consider the following simplifications for the filter
parameters:
\begin{align*}
u_t^i & := u_t\\
\sum_{j=1}^N(b_T^{ij})^\top X_0^j & := 
c_T^\top (X_0^i -  X_0^{(N)}) + b_T^\top X_0^{(N)} 
\end{align*}
where $X_0^{(N)}:=\frac{1}{N} \sum_{i=1}^N X_0^i$ is the empirical mean.  Using this new parametrization, the filter~\eqref{eq:FPF_LP_est_1} becomes
\begin{equation}
S_T^i = 
c_T^\top (X_0^i -  X_0^{(N)}) + b_T^\top X_0^{(N)} + \int_0^T u_t^\top \ud Z_t
\label{eq:FPF_LP_est_2}
\end{equation}

\newP{Step 2} The mean-field limit of the estimator is obtained by
letting $N\rightarrow \infty$.  In the limit, the random variable is
denoted as $\bar{S}_T$.  The mean-field counterpart of the
filter~\eqref{eq:FPF_LP_est_2}  is
\begin{equation}
\bar{S}_T =
c_T^\top (\bar{X}_0- m_0) +   b_T^\top m_0 + \int_0^T u_t^\top \ud Z_t
\label{eq:FPF_LP_est}
\end{equation}
The optimization problem is to chose the filter parameters $b_T$,
$c_T$ and $u_t$ to minimize the following mean-squared cost
\begin{equation*}
\min_{b_T,c_T,u_t} \;\; \E [|\E[f(\bar{S}_T)|{\cal Z}_T] - f(a^\top X_T)|^2]
\label{eq:multi_objective}
\end{equation*}

\newP{Step 3} In the final step, we restrict the class of functions to
$f(x)=x$ and $f(x)=x^2$. The justification for considering only linear
and quadratic functions is that the density of the random variable
$\bar{S}_T$ is Gaussian (see also the
Remark~\ref{remark:aposteriori_justification} at the end of this section).  

\subsection{Optimization problem}

In summary, the mathematical problem is a bi-objective optimization
problem:
\begin{equation}\label{eq:MOO_function}
\min_{b_T,c_T,u_t} \left( \E [|\E[\bar{S}_T|{\cal Z}_T] - a^\top X_T|^2],\;
\E [|\E[\bar{S}_T^{\,2}|{\cal Z}_T] - (a^\top X_T)^2|^2]\right)
\end{equation}
subject to the respective constraints for the state, observation, and the estimator:
\begin{subequations}\label{eq:opt_system_def}	
\begin{align}
\ud X_t &=  A_tX_t \ud t + \ud B_t,\quad X_0\sim {\cal N}(m_0,\Sigma_0)\label{eq:dyn*}\\
\ud Z_t &= C_tX_t\ud t + \ud W_t\label{eq:obs*} \\
\bar{S}_T & = b_T^\top m_0+
c_T^\top (\bar{X}_0- m_0) +  \int_0^T u_t^\top \ud Z_t, \quad \bar{X}_0\sim {\cal N}(m_0,\Sigma_0)\label{eq:opt_system_def_c}
\end{align} 
\end{subequations}
where the noise processes $W_t$ and $B_t$, and the initial conditions $X_0$ and
$\bar{X}_0$, are all assumed to be mutually independent; recall that $\clZ_t:=\sigma(Z_s;s\in[0,t])$
denotes the time-history of observations up to time $t$ (filtration).   

A solution to this problem appears in the following section.  

\section{Deterministic particle filter via duality}
\label{sec:det}

Denote $\Sigma_t$ to be the solution of the DRE~\eqref{eq:KF-variance}
with the initial condition $\Sigma_0$, and $\k_t:=\Sigma_tC^\top
R_t^{-1}$ is the Kalman gain.  

Define the following state transition matrices:
\begin{equation}\label{eq:Phi}
\dfrac{\ud \Phi}{\ud t}(t;\tau) = \big(-A_t^\top + C_t^\top \k_t^\top \big)\Phi(t;\tau)
\end{equation}
\begin{equation}\label{eq:Psi}
\dfrac{\ud\Psi}{\ud t}(t;\tau) = \big(-A_t^\top + \dfrac{1}{2}C_t^\top
\k_t^\top -\dfrac{1}{2}\Sigma_t^{-1}Q_t\big)\Psi(t;\tau)
\end{equation}
with $\Phi(\tau;\tau)=\Psi(\tau;\tau)=I$, the identity matrix. 

These definitions are useful to describe the solution to the optimization
problem~\eqref{eq:MOO_function}, as presented in the following
Theorem.  Its proof appears in the Appendix (\Sec{apdx:pf_thm1}).   

\medskip

\begin{theorem}\label{thm:thm1}
Consider the optimization problem~\eqref{eq:MOO_function} subject to
the constraints~\eqref{eq:opt_system_def}. Its solution is as
follows:
\begin{subequations}\label{eq:opt_solution}
\begin{align}
b_T &= \Phi(0;T)a\\
c_T &= \Psi(0;T)a\\
u_t &= \k_t^\top\Phi(t;T)a
\end{align}
\end{subequations}
The solution for $c_T$ is unique up to a sign for the scalar ($d=1$)
case but it is not unique for the vector ($d>1$) case.  

Using the optimal parameters~\eqref{eq:opt_solution}, the estimator is 
\begin{equation}\label{eq:opt_S_T}
\begin{split}
\bar{S}_T =& a^\top\Big(\Phi^\top(0;T)m_0 + \Psi^\top(0;T)(\bar{X}_0-m_0)\\&+\int_0^T \Phi^\top(t;T)\k_t\ud Z_t\Big),\quad \bar{X}_0\sim {\cal N}(m_0,\Sigma_0)
\end{split}
\end{equation}
This estimator is exact.  That is,
\[
{\sf E}(g(\bar{S}_T)|\clZ_T) = {\sf E}(g(a^\top X_T)|\clZ_T) 
\]
for all smooth test functions $g$ and vectors $a\in\Re^d$. 
\end{theorem}

\medskip

The following Proposition provides the differential (recursive) form
of the filter which also reveals the connection to the FPF algorithm
(compare with formula~\eqref{eq:Xbart} for the deterministic FPF). The proof appears in the Appendix (\Sec{apdx:pf_thm2}).

\medskip

\begin{proposition}\label{thm:thm2}
For any given $a\in\Re^d$ and $T\geq 0$, $\bar{S}_T = a^\top
\bar{X}_T$ where $\bar{X}_T$ is the strong solution of the following
mean-field sde:
\begin{align}
\ud &\bar{X}_t = A_t \bar{X}_t\ud t + \k_t\Big(\ud Z_t - \frac{C_t\bar{X}_t + C_t\bar{m}_t}{2}\ud t\Big) + \frac{1}{2}Q_t \Sigma_t^{-1}(\bar{X}_t-\bar{m}_t) \ud t \nonumber
\end{align}
where $\bar{m}_t = \E[\bar{X}_t|\clZ_t]$ and the initial condition
$\bar{X}_0\sim{\cal N}(m_0,\Sigma_0)$.  
\end{proposition}

 \medskip

\begin{remark}\label{remark:aposteriori_justification}
The derivation of the deterministic FPF offers an {\em a posteriori}
justification of our simplifying choices for the filter parameters (in
step 1) and the functions $f$ (in step 3).  That these choices were
sufficient is primarily due to the linear Gaussian nature of the
problem.  More generally, if the distribution is sub-Gaussian, the
method of moments
suggests considering $\{f(x) = x^k \; |\; k = 1,2,\ldots\}$ as a class of
functions~\cite[Section 30]{billingsley1986probability}.  
\end{remark}



\section{Stochastic particle filter via duality}
\label{sec:stoch}

In this section, we extend the deterministic filter
structure~\eqref{eq:FPF_LP_est_2} to now include noise terms:
\begin{equation*}\label{eq:stoc_filter}
\begin{split}
{S}_T^i =& b_T^\top X_0^{(N)}+
c_T^\top (X_0^i- X_0^{(N)})\\ &+ \int_0^T u_t^\top \ud Z_t  + \int_0^T v_t^\top \ud B_t^i + \int_0^T w_t^\top \ud W_t^i
\end{split}
\end{equation*}
where the $\{B_t^i\}$ and $\{W_t^i\}$ are independent copies of the
process noise (with covariance $Q_t$) and the measurement error (with
covariance $R_t$), respectively.  These are also independent of
$B_t,W_t,X_0,X_0^i$.  The design problem is to chose the filter
parameters $\{b_T,c_T,u_t,v_t,w_t\}$.

The motivation of considering this type of
filter structure is two-fold: (i) ${S}_T^i$ now depends linearly on
{\em all} of the data---observation, initial condition as well as the
copies of the process noise and the measurement noise; and (ii) 
particle filters with noise terms are widely used in practice~\cite{TaghvaeiASME2017,bergemann2012}.  For
example, the stochastic FPF~\eqref{eq:Xbart-stochastic} includes noise term to simulate the effect
of process noise.  The ensemble Kalman filter 
include noise terms to simulate the effect of both process noise and measurement noise~\cite{bergemann2012}.

On repeating the steps~1-3, the stochastic counterpart of the
optimization problem~\eqref{eq:MOO_function}-\eqref{eq:opt_system_def}
is obtained as follows:
\begin{equation}\label{eq:MOO_2}
\min_{\substack{b_T,c_T,\\u_t,v_t,w_t}} \left( \E [|\E[\bar{S}_T|{\cal Z}_T] - a^\top X_T|^2],\;
\E [|\E[\bar{S}_T^{\,2}|{\cal Z}_T] - (a^\top X_T)^2|^2]\right)
\end{equation}
subject to
\begin{subequations}\label{eq:opt_2}
\begin{align}
\ud X_t =&  A_tX_t \ud t + \ud B_t,\quad X_0\sim {\cal N}(m_0,\Sigma_0)\\
\ud Z_t =& C_tX_t\ud t + \ud W_t\\
\bar{S}_T =& b_T^\top m_0+
c_T^\top (\bar{X}_0- m_0) + \int_0^T u_t^\top \ud Z_t \nonumber\\ & + \int_0^T v_t^\top \ud \bar{B}_t + \int_0^T w_t^\top \ud \bar{W}_t,\quad \bar{X}_0\sim {\cal N}(m_0,\Sigma_0)\label{eq:opt_2_S}
\end{align} 
\end{subequations}

\medskip

A solution of the optimization problem is given in the following
theorem whose proof appears in the Appendix (\Sec{apdx:pf_thm3}).  
As in the deterministic case, the solution is given in terms of certain state transition matrices.  The
state-transition matrix $\Phi$ is as defined in~\eqref{eq:Phi}.  The state
transition matrix $\Psi$ is modified to the following:
\begin{equation*}
\dfrac{\ud\Psi}{\ud t}(t;\tau) = \Big(-A_t^\top  +
\frac{1+\gamma_2^2}{2}C_t^\top \k_t^\top-
\frac{1-\gamma_1^2}{2}\Sigma_t^{-1}Q_t \Big) \Psi(t;\tau)
\end{equation*}  
where $\gamma_1$ and $\gamma_2$ are real-valued parameters.  With
$\gamma_1=\gamma_2=0$, one obtains the original
definition~\eqref{eq:Psi} for $\Psi$. 

\medskip

\begin{theorem}\label{thm:thm3}
Consider the optimization problem~\eqref{eq:MOO_2} subject to the
dynamic constraints~\eqref{eq:opt_2}.  One solution of the
optimization problem is as follows:
\begin{align*}
	b_T &= \Phi(0;T)a\\
	c_T &= \Psi(0;T)a\\
	u_t &= \k_t^\top\Phi(t;T)a\\
	v_t &= \gamma_1\Psi(t;T)a\\
	w_t &= \gamma_2\k_t^\top\Psi(t;T)a
\end{align*}
where the values of the parameters $\gamma_1$ and $\gamma_2$ can be
arbitrarily chosen.  The resulting optimal estimator is:
	\begin{equation}
	\begin{split}
	\bar{S}_T =& a^\top\Big(\Phi^\top(0;T)m_0+
        \Psi^\top(0;T)(\bar{X}_0-m_0)\\&+\int_0^T\Phi^\top(t;T)\k_t\ud
        Z_t + \gamma_1\int_0^T \Psi^\top(t;T)\ud \bar{B}_t \\&+ \gamma_2\int_0^T\Psi^\top(t;T)\k_t\ud \bar{W}_t\Big)
	\end{split}
	\end{equation}
The estimator is exact for every choice of the parameter values $\gamma_1$ and $\gamma_2$.  That is,
\[
{\sf E}(g(\bar{S}_T)|\clZ_T) = {\sf E}(g(a^\top X_T)|\clZ_T) 
\]
for all smooth test functions $g$ and vectors $a\in\Re^d$. 
\end{theorem}

\medskip

The differential form of the filter is given in the following
Proposition whose proof appears in the Appendix.   

\medskip

\begin{proposition}
For any given $a\in\Re^n$ and $T\geq 0$, $\bar{S}_T = a^\top
\bar{X}_T$ where $\bar{X}_T$ is the strong solution of the following
mean-field sde:
\begin{equation*}
\begin{split}
\ud \bar{X}_t =& A_t \bar{X}_t\ud t+ \gamma_1\ud \bar{B}_t + \frac{1-\gamma_1^2}{2}Q_t \Sigma_t^{-1}(\bar{X}_t-\bar{m}_t) \ud t \\&+ \k_t\Big(\ud Z_t - C_t\Big(\frac{(1+\gamma_2^2)\bar{X}_t + (1-\gamma_2^2)\bar{m}_t}{2}\Big)\ud t+\gamma_2\ud\bar{W}_t\Big)
\end{split}
\end{equation*}
where $\bar{m}_t = \E[\bar{X}_t|\clZ_t]$ and the initial condition
$\bar{X}_0\sim{\cal N}(m_0,\Sigma_0)$.  
\end{proposition}

\medskip

\begin{remark}
The parameters $\gamma_1$ and $\gamma_2$ parametrize a homotopy of
ensemble Kalman filters/linear FPFs, all of which are exact in the
linear Gaussian settings.
\begin{enumerate}
\item For $\gamma_1=\gamma_2=0$, one obtains the deterministic form of
  the linear FPF (Eq.~\eqref{eq:Xbart}).
\item For $\gamma_1=1$ and $\gamma_2=0$, one obtains the stochastic
  linear FPF (Eq.~\eqref{eq:Xbart-stochastic}).  
\item For $\gamma_1=\gamma_2=1$, one obtains the original form of the
  ensemble Kalman filter~\cite[Eq.~(7)]{TaghvaeiASME2017} where a copy of measurement noise is
  introduced in the error.   
\end{enumerate}
By varying the parameters $\gamma_1$ and $\gamma_2$ in the range
$[0,1]$, one goes from the two stochastic filters to the deterministic
filter.  The filters are exact for
arbitrary  (even time-varying) values of parameters $\gamma_1$ and
$\gamma_2$.

\end{remark}

\section{Conclusion and directions for future work}

In this paper, a novel homotopy of exact linear Gaussian particle filters
is derived, based on a certain extension of the classical duality.
The filter is also related to the linear FPF 
and the ensemble Kalman filter.

There are several possible directions of future work. The first
direction is to extend the proposed duality framework to now
incorporate the finite-$N$ effects. Although it was not an explicit
focus of this paper, the relationship to optimal control potentially
provides new tools to analyze and interpret the optimality properties
of the particle filter, in terms of the simulation variance of its
estimates.  The other direction is to extend
the duality framework to nonlinear and non-Gaussian settings. Note
that FPF algorithm is known to be exact in these settings. If
possible, it will certainly be of interest to derive the general FPF
using the duality framework.


\bibliographystyle{IEEEtran}
\bibliography{ref2}

\appendix



\subsection{Proof of Theorem 1}\label{apdx:pf_thm1}

Express the filter~\eqref{eq:FPF_LP_est} as
\[
\bar{S}_T =
c_T^\top (\bar{X}_0- m_0) +  M_T(b_T,u_t)
\]
where 
\begin{equation}\label{eq:hat_S}
M_T(b_T,u_t):= b_T^\top m_0 + \int_0^T u_t \ud Z_t = \E[\bar{S}_T|{\cal Z}_T]
\end{equation}
is the conditional mean of $\bar{S}_T$.  

The first objective involves minimization of the mean-squared error
(filter variance):
\[
\E [|\E[\bar{S}_T|{\cal Z}_T] - a^\top X_T|^2] \;\; = \;\; \underbrace{\E [|M_T - a^\top X_T|^2]}_{\text{term (i)}}
\]

The second objective function in~\eqref{eq:MOO_function} is
expressed as a sum of two terms:
\begin{align}
\E\Big[\big|\E\big[&\bar{S}_T^2|\clZ_T\big]-(a^\top X_T)^2\big|^2\Big]
= \E\Big[\big|c_T^\top\Sigma_0c_T + M_T^2 - (a^\top X_T)^2\big|^2\Big]\nonumber\\
=&\;\;\underbrace{{\sf var}\big[(a^\top X_T)^2-M_T^2\big]}_{\text{term (ii)}} \;\;+\;\; \underbrace{\big(c_T^\top\Sigma_0c_T - \E\big[(a^\top
  X_T)^2-M_T^2\big]\big)^2}_{\text{term (iii)}}\label{eq:op_term23}
\end{align}

The optimization of the three terms~(i)-(iii) is the subject of the
three steps in this proof.  
\begin{enumerate}
\item In step 1, term~(i) is minimized by choosing
  the vector $b_T$ and the control $u_t$.  The solution of this
  problem is given by classical duality.
\item In step 2, it is shown that solution thus obtained also
  minimizes term~(ii).  
\item In step 3, the vector $c_T$ is chosen to minimize the square
  term~(iii).  Its minimum value is $0$.    
\end{enumerate}
The details of the three steps appear next.



\newP{Step 1} The solution to the term~(i) minimization problem is
given by classic duality~\cite{astrom1970}.  A dual process is
introduced: 
\begin{equation}\label{eq:dual_pf}
\ud y_t = -A_t^\top y_t \ud t + C_t^\top u_t \ud t,\quad y_T = a
\end{equation}
Since $\ud(y_t^\top X_t) = u_t^\top C_tX_t \ud t + y_t^\top \ud B_t $,
$$
a^\top X_T = y_0^\top X_0 + \int_0^T u_t^\top C_tX_t\ud t + \int_0^T y_t^\top \ud B_t
$$
Using~\eqref{eq:hat_S} and $\ud Z_t = C_tX_t \ud t + \ud W_t$, one obtains
$$
M_T-a^\top X_T = b_T^\top m_0 - y_0^\top X_0 + \int_0^T u_t^\top \ud W_t - \int_0^T y_t^\top \ud B_t
$$
Squaring and taking expectations, the term~(i) is expressed as
\begin{equation*}
\big((b_T^\top - y_0^\top)m_0\big)^2+ y_0^\top \Sigma_0 y_0 + \int_0^T y_t^\top Q_t y_t + u_t^\top R_t u_t \ud t
\end{equation*}

With the minimizing choice of $b_T = y_0 =: b_T^*$, the term~(i) minimization
problem is transformed into a linear quadratic (LQ) optimal control problem 
\[
\min_{u_t}\quad y_0^\top \Sigma_0 y_0 + \int_0^T y_t^\top Q_t y_t + u_t^\top R_t u_t \ud t
\]
subject to the dynamic constraints given by the dual
system~\eqref{eq:dual_pf}.  

The optimal control law is easily obtained in the feedback form 
\[
u_t = \k_t^\top y_t 
\]
where the gain $\k_t:=\Sigma_t C^\top
R_t^{-1}$ is obtained by solving the DRE~\eqref{eq:KF-variance}.  
Upon using the defintion~\eqref{eq:Phi} of the state transition matrix
$\Phi(t,\tau)$, 
\[
u_t = \k_t^\top y_t = \Phi(t;T)y_T = \Phi(t;T)a =: u_T^*
\]

We use the notation $\hat{S}_T$ to denote $M_T$ with the optimal
choice of the parameters $b_T^*$ and $u_t^*$.  Explicitly,
\begin{equation}\label{eq:hatS}
\hat{S}_T := M_T(b_T^*,u_t^*) = a^\top\Phi^\top(0;T)m_0 + \int_0^Ta^\top\Phi^\top(t;T)\k_t\ud Z_t
\end{equation}

$\hat{S}_T$ has the following properties which
are useful in the remainder of the proof:
\begin{enumerate}
\item $\hat{S}_T$  is the conditional mean, i.e., $\hat{S}_T = E[a^\top
  X_T|{\cal Z}_T]$.  
\item The error $(a^\top X_T-\hat{S}_T)$ satisfies the orthogonal
  property whereby
$\E\big[(a^\top X_T-\hat{S}_T)\hat{S}_T\big] = 0$.  
\item The optimal value $\E\big[|\hat{S}_T-a^\top X_T|^2\big] = a^\top
  \Sigma_T a$.
\end{enumerate}

\newP{Step 2} The optimal parameters, $b_T^*$ and $u_t^*$, obtained in step~1 also 
minimize the term~(ii). 
The reasoning is as follows:
Note that the term~(ii) involves minimization of the variance again by choosing
$b_T$ and $u_t$.  Since the resulting random variable $M_T^2$ is   
${\cal Z}_T$-measurable, any minimizer of
term~(ii) is of the general form
\begin{equation}\label{eq:term_2_opt}
M_T^2 = \E[(a^\top X_T)^2|\clZ_T] + \text{(const.)}
\end{equation}
Now, the conditional variance of $a^\top X_T$ is given by
\begin{align*}
{\sf var}[a^\top X_T|\clZ_T] &=\E[(a^\top X_T-\E[a^\top X_T|\clZ_T])^2|\clZ_T]\\
&=\E[(a^\top X_T-\E[a^\top X_T|\clZ_T])^2]\\
&=\E[(a^\top X_T-\hat{S}_T)^2] = a^\top \Sigma_T a
\end{align*}
The second equality comes from the fact that orthogonality implies
independence for a Gaussian random variable. Since
\begin{align*}
{\sf var}[a^\top X_T|\clZ_T] &=\E[(a^\top X_T)^2|\clZ_T] - (\E[a^\top X_T|\clZ_T])^2\\
&=\E[(a^\top X_T)^2|\clZ_T] - \hat{S}_T^2
\end{align*}
we have
$$
\hat{S}_T^2 = \E[(a^\top X_T)^2|\clZ_T] - a^\top\Sigma_T a
$$
Comparing with~\eqref{eq:term_2_opt}, it follows that $M_T=\hat{S}_T$
is a minimizer of~the term~(ii).  Consequently, $b_T^*$ and $u_t^*$
are optimal choices for minimizing term~(ii) as well.  It is noted
that $\hat{S}_T$ (defined in~\eqref{eq:hatS}) is the {\em only} solution that simultaneously 
minimizes both terms~(i) and~(ii).

%
%
%

\newP{Step 3} Now that the vector $b_T$ and the control $u_t$ have
been obtained, the vector $c_T$ is chosen to make the square term~(iii) zero.  We have 
\begin{align*}
c_T^\top\Sigma_0c_T &= \E\big[(a^\top X_T)^2-(\hat{S}_T)^2\big] \\
&=\E\big[(a^\top X_T)^2-(\hat{S}_T)^2\big] - 2 \E\big[(a^\top X_T-\hat{S}_T)\hat{S}_T\big] \\
&= \E\big[|a^\top X_T - \hat{S}_T|^2\big]  = a^\top \Sigma_T a 
\end{align*}
Therefore, the minimizing choice of $c_T$ is obtained by solving the
scalar equation
\[
c_T^\top\Sigma_0c_T = a^\top \Sigma_T a 
\]
It is straightforward to verify that $c_T = \Psi(0;T)a$ is a solution
of this equation.  A more constructive proof follows from introducing
a backward-time process
\[
\frac{\ud \xi_t}{\ud t} = \bigg(-A_t^\top +\frac{1}{2}C_t^\top
\k_t^\top-\frac{1}{2}\Sigma_t^{-1}Q_t\bigg) \xi_t, \quad \xi_T=a
\] 
Since $\Sigma_t$ is a solution of the DRE~\eqref{eq:KF-variance}, it is easy to then verify that $\ud(\xi_t^\top \Sigma_t \xi_t)
= 0$, and so $\xi_0^\top \Sigma_0 \xi_0 = a^\top \Sigma_T a$.  Thus,
$c_T = \xi_0= \Psi(0;T)a$ is a solution.

The proof for exactness is deferred to the following Sections.\qed



\subsection{Proof of Proposition 2}\label{apdx:pf_thm2}

The equation~\eqref{eq:opt_S_T} for the optimal estimator has the
following form: $\bar{S}_T = a^\top \bar{X}_T$.  Upon denoting the
(arbitrary) final time $T$ simply as $t$, one
writes
\begin{equation}\label{eq:X_t}
\bar{X}_t =\Psi^\top(0;t)(\bar{X}_0-m_0) + \bar{m}_t
\end{equation}
where $\bar{m}_t = \Phi^\top(0;t)m_0 +\int_0^t \Phi^\top(s;t)\k_s\ud
Z_s =\E[\bar{X}_t|{\cal Z}_t] $.

Upon differentiating~\eqref{eq:X_t}, and using the
formulae~\eqref{eq:Phi} and~\eqref{eq:Psi} for the state transition matrices,
\begin{align*}
\ud \bar{X}_t =& (A_t-\k_tC_t)\Phi^\top(0;t)m_0 \ud t\\&+ \Big(A_t-\frac{1}{2}\k_tC_t+\frac{1}{2}Q_t\Sigma_t^{-1}\Big)\Psi^\top(0;t)(\bar{X}_0-m_0) \ud t\\&+\Big(\int_0^t (A_t-\k_tC_t)\Phi^\top(s;t)\k_s\ud Z_s\Big) \ud t + \k_t \ud Z_t\\
=&(A_t-\k_tC_t)\bar{X}_t \ud t + \Big(\frac{1}{2}\k_tC_t+\frac{1}{2}Q_t\Sigma_t^{-1}\Big)(\bar{X}_t-\bar{m}_t) \ud t
\end{align*}
This yields the recursive formula for the filter.  The proof of
exactness of this filter has already appeared in~\cite{AmirACC2016} (see also~\Prop{thm:consistency-FPF-lin}).  It is also
a special case of the more general stochastic filter whose exactness proof
appears in the \Sec{apdx:pf_exact}.    

\medskip

\subsection{Proof of Theorem 2 and Proposition 3}\label{apdx:pf_thm3}
Express the filter~\eqref{eq:opt_2_S} as
$$
\bar{S}_T = M_T + c_T^\top (\bar{X}_0-m_0)+\int_0^T v_t^\top \ud \bar{B}_t + \int_0^T w_t^\top \ud \bar{W}_t
$$
where $M_T= M_T(b_T,u_t) =  \E[\bar{S}_T|{\cal Z}_T]$ is as defined
in~\eqref{eq:hat_S}.  The other three terms model
the effects of randomness due to the initial condition, process noise
and the measurement noise, respectively.

Since both the noise terms have zero-mean, these do not affect the
terms~(i) and~(ii) as introduced in the proof of
\Theorem{thm:thm1}.  The counterpart of~\eqref{eq:op_term23} now is:
\begin{align*}
& \E\Big[\big|\E\big[\bar{S}_T^2|\clZ_T\big]-(a^\top X_T)^2\big|^2\Big]
=\;\underbrace{{\sf var}\big[(a^\top
  X_T)^2-M_T^2\big]}_{\text{term (ii)}}\; \\ & 
+\; \underbrace{\big(c_T^\top\Sigma_0c_T +\int_0^T v_t^\top Q_tv_t +  w_t^\top R_t w_t \ud t - \E\big[(a^\top
	X_T)^2-M_T^2\big]\big)^2}_{\text{term (iii)}'}
\end{align*}

Since the terms~(i) and~(ii) are identical to the ones in the proof the \Theorem{thm:thm1}, the
steps 1 and 2 apply in an identical manner.  The optimal solution for $b_T$
and $u_t$ is thus the same as before.  The optimal $M_T$
is given by $\hat{S}_T$ (see~\eqref{eq:hatS}).

Upon setting the square term~(iii)$'$ to its minimum value, zero, the scalar equation for $c_T$ is now given by
\begin{equation}\label{eq:opt_cond_ct}
c_T^\top\Sigma_0c_T = a^\top \Sigma_T a-\int_0^Tv_t^\top Q_t v_t + w_t^\top R_t w_t \ud t
\end{equation}
This equation has many solutions.  We pick one solution by introducing
the following modification of the process $\xi_t$:
\begin{equation*}\label{eq:dyn_xi_2}
\frac{\ud \xi_t}{\ud t} = \left(-A_t^\top  + \frac{1+\gamma_2^2}{2}C_t^\top \k_t^\top- \frac{1-\gamma_1^2}{2}\Sigma_t^{-1}Q_t\right) \xi_t,\; \xi_T = a
\end{equation*}
where $\gamma_1$ and $\gamma_2$ are arbitrary constants. Then 
\begin{equation*}
\ud (\xi_t^\top \Sigma_t \xi_t) = \gamma_1^2\xi_t^\top Q_t\xi_t \ud t + \gamma_2^2\xi_t^\top \k_t R_t \k_t^\top \xi_t\ud t,\quad \xi_T = a
\end{equation*}
Therefore, upon setting $v_t = \gamma_1\xi_t$ and $w_t =
\gamma_2\k_t^\top\xi_t$, one obtains
$$
\xi_0^\top \Sigma_0 \xi_0 = a^\top \Sigma_T a-\int_0^Tv_t^\top Q_t v_t + w_t^\top R_t w_t \ud t
$$
and thus
\begin{align*}
c_T &= \xi_0 = \Psi(0;T)a\\
v_t &= \gamma_1 \Psi(t;T)a\\
w_t &= \gamma_2 \k_t^\top \Psi(t;T)a
\end{align*}
solves~\eqref{eq:opt_cond_ct}.

Using these parameters, the estimator is
$\bar{S}_T = a^\top \bar{X}_T$
where
\begin{equation*}
\begin{split}
\bar{X}_t = & \Phi^\top(0;t)m_0+ \Psi^\top(0;t)(\bar{X}_0-m_0)\\&+\int_0^t\Phi^\top(s;t)\k_s\ud Z_s + \gamma_1\int_0^t \Psi^\top(s;t)\ud \bar{B}_s \\&+ \gamma_2\int_0^t\Psi^\top(s;t)\k_s\ud \bar{W}_s
\end{split}
\end{equation*}
Its differential form is easily obtained as
\begin{equation}\label{eq:gen_fpf}
\begin{split}
\ud \bar{X}_t =& A_t \bar{X}_t\ud t+ \gamma_1\ud \bar{B}_t + \frac{1-\gamma_1^2}{2}Q_t \Sigma_t^{-1}(\bar{X}_t-\bar{m}_t) \ud t \\&+ \k_t\Big(\ud Z_t - C_t\big(\frac{(1+\gamma_2^2)\bar{X}_t + (1-\gamma_2^2)\bar{m}_t}{2}\big)\ud t+\gamma_2\ud\bar{W}_t\Big)
\end{split}
\end{equation}

\subsection{Proof of the exactness of the filter}\label{apdx:pf_exact}

The filter~\eqref{eq:gen_fpf} is a linear sde, the initial condition $\bar{X}_0$
is Gaussian, and the noise terms are also Gaussian.  Therefore, the
solution $\bar{X}_t$ is Gaussian for all $t > 0$.  Thus, in order to show
exactness, all we need to show is that the equations for the
conditional mean and the variance evolve according to the Kalman
filter equations.     

Upon taking a conditional expectation of~\eqref{eq:gen_fpf},
\begin{equation*}\label{eq:hat_evolve}
\ud \bar{m}_t = A_t\bar{m}_t + \k_t(\ud Z_t - C_t\bar{m}_t\ud t),\quad \bar{m}_0 = m_0
\end{equation*}
This is the same as the equation~\eqref{eq:KF-mean} of the Kalman
filter. 

Next define the error process $e_t:= \bar{X}_t - \bar{m}_t$.  The
equation for the error process is given by
\begin{equation*}
\begin{split}
\ud e_t =& A_t e_t\ud t + \frac{1-\gamma_1^2}{2}Q_t \Sigma_t^{-1}e_t \ud t \\&
- \k_tC_t\frac{1+\gamma_2^2}{2}e_t\ud t+ \gamma_1\ud \bar{B}_t +\gamma_2\k_t\ud\bar{W}_t
\end{split}
\end{equation*}

Upon squaring the taking expectations, one obtains the differential
equation for the variance $\bar{\Sigma}_t = {\sf var}(\bar{X}_t)$ as follows:
\begin{equation*}
\begin{split}
\frac{\ud \bar{\Sigma}_t}{\ud t} =& A_t \bar{\Sigma}_t + \bar{\Sigma}_t A_t^\top  + \frac{1-\gamma_1^2}{2}(Q_t \Sigma_t^{-1}\bar{\Sigma}_t +\bar{\Sigma}_t\Sigma_t^{-1} Q_t) \\&- \frac{1+\gamma_2^2}{2}(\k_tC_t\bar{\Sigma}_t-\bar{\Sigma}_tC^\top \k_t^\top)+ \gamma_1^2 Q_t +\gamma_2^2\k_t R_t\k_t^\top 
\end{split}
\end{equation*}
with initial condition $\bar{\Sigma}_0 = \Sigma_0$.  
This is a linear ode (in $\bar{\Sigma}_t$).  It thus admits a unique
solution.  It is straightforward to verify that $\bar{\Sigma}_t =
\Sigma_t$ in fact solves the equation (recall that $\Sigma_t$ is the
solution to the DRE~\eqref{eq:KF-variance} with initial condition $\Sigma_0$).  

Therefore, the
filter~\eqref{eq:gen_fpf} is exact for all time $t \geq 0$.


\qed

\end{document}